\documentclass[11pt]{amsart}
\usepackage{amssymb, amstext, amscd, amsmath}
\usepackage{mathtools, xypic, color, dsfont, rotating}
\usepackage{verbatim}
\usepackage{enumitem}

\usepackage{showlabels}

\usepackage{hyperref}

\makeatletter
\def\@cite#1#2{{\m@th\upshape\bfseries%
[{#1\if@tempswa{\m@th\upshape\mdseries, #2}\fi}]}}
\makeatother
\theoremstyle{plain}
\newtheorem{theorem}{Theorem}[section]
\newtheorem{corollary}[theorem]{Corollary}
\newtheorem{proposition}[theorem]{Proposition}

\theoremstyle{definition}
\newtheorem{definition}[theorem]{Definition}
\newtheorem{example}[theorem]{Example}

\theoremstyle{remark}

\mathtoolsset{centercolon}
  \newcommand{\B}{{\mathcal{B}}}

  \newcommand{\F}{{\mathcal{F}}}

  \newcommand{\I}{{\mathcal{I}}}

  \newcommand{\M}{{\mathcal{M}}}
  
\renewcommand{\O}{{\mathcal{O}}}

  \newcommand{\R}{{\mathcal{R}}}

  \newcommand{\T}{{\mathcal{T}}}
  \newcommand{\U}{{\mathcal{U}}}

\renewcommand{\phi}{\varphi}
\def\si{\sigma}
\def\al{\alpha}
\def\be{\beta}

\def\ga{\gamma}
\newcommand\vpi{\varphi}

\newcommand{\rC}{{\mathrm{C}}}

\newcommand{\bD}{\mathbb{D}}
\newcommand{\bF}{\mathbb{F}}

\newcommand{\bT}{\mathbb{T}}
\newcommand{\bZ}{\mathbb{Z}}
\newcommand{\bR}{\mathbb{R}}
\newcommand{\fA}{{\mathfrak{A}}}

\newcommand{\Be}{{\mathbf{e}}}

\newcommand{\Bi}{{\mathbf{i}}}

\newcommand{\Bn}{{\mathbf{n}}}

\newcommand{\Bs}{{\mathbf{s}}}
\newcommand{\Bt}{{\mathbf{t}}}

\newcommand{\AND}{\text{ and }}
\newcommand{\FOR}{\text{ for }}
\newcommand{\FORAL}{\text{ for all }}

\newcommand{\qand}{\quad\text{and}\quad}

\newcommand{\qfor}{\quad\text{for}\ }
\newcommand{\qforal}{\quad\text{for all}\ }

\newcommand{\ca}{\mathrm{C}^*}

\newcommand{\cenv}{\mathrm{C}^*_{\textup{env}}}

\newcommand{\lip}{\langle}
\newcommand{\rip}{\rangle}
\newcommand{\ip}[1]{\lip #1 \rip}

\newcommand{\ol}{\overline}
\newcommand{\wt}{\widetilde}
\newcommand{\wh}{\widehat}

\newcommand{\ad}{\operatorname{ad}}

\newcommand{\Aut}{\operatorname{Aut}}

\newcommand{\dirlim}{\varinjlim}
\newcommand{\diag}{\operatorname{diag}}
\newcommand{\Dim}{\operatorname{dim}}

\newcommand{\End}{\operatorname{End}}

\newcommand{\id}{{\operatorname{id}}}

\newcommand{\iso}{{\operatorname{is}}}

\newcommand{\mt}{\emptyset}
\newcommand{\nc}{\operatorname{nc}}
\newcommand{\nd}{\operatorname{nd}}

\newcommand{\supp}{\operatorname{supp}}
\newcommand{\uni}{{\operatorname{un}}}
\newcommand{\bo}[1]{\mathbf{#1}} 
\newcommand{\un}[1]{{\underline{#1}}} 

\addtocontents{toc}{\protect\setcounter{tocdepth}{1}}

\begin{document}


\title[Semicrossed Products of Operator Algebras]{Semicrossed Products of Operator Algebras:\\ A Survey}
\author[K.R. Davidson]{Kenneth R. Davidson}
\address{Pure Mathematics Department\\University of Waterloo\\Waterloo, ON\\Canada \ N2L--3G1}
\email{krdavids@uwaterloo.ca}
\thanks{First author partially supported by NSERC, Canada.}
\author[A.H. Fuller]{Adam H. Fuller}
\address{Mathematics Department\\University of Nebraska-Lincoln\\Lincoln, NE\; 68588--0130\\USA}
\email{afuller7@math.unl.edu}

\author[E.T.A. Kakariadis]{Evgenios T.A. Kakariadis}
\address{Department of Mathematics\\Ben-Gurion University of the Negev\\Be'er Sheva 84105\\Israel}
\email{kakariad@math.bgu.ac.il}
\thanks{Third author partially supported by the Kreitman Foundation Post-doctoral Fellow Scholarship, and the Skirball Foundation via the Center for Advanced Studies in Mathematics at Ben-Gurion University of the Negev.}


\begin{abstract}
Semicrossed product algebras have been used to study dynamical systems since their introduction by Arveson in 1967. In this survey article, we discuss the history and some recent work, focussing on the conjugacy problem, dilation theory and C*-envelopes, and some connections back to the dynamics.
\end{abstract}

\dedicatory{We dedicate this article to the memory of William B. Arveson,\\ who has inspired all of us with his vision and deep insights\\ that completely changed the way we look at operator theory.}%
\maketitle

\section{Introduction}

The study of operator algebras arising from dynamical systems is almost as old as the study of operator algebras themselves. W*-crossed products were originally studied in the seminal work of Murray and von Neumann. An analogous theory of C*-crossed products has also been developed to build a universal C*-algebra which encodes the action of a \emph{group} of $*$-automorphisms on a C*-algebra $A$. See Pedersen \cite{Ped79} or Williams \cite{Wil07} for an introduction. We instead turn our focus to actions of \emph{semigroups} on an \emph{arbitrary} operator algebra $A$ by \emph{endomorphisms}. In this case, the natural object encoding this action is a nonselfadjoint operator algebra, even when $A$ is a C*-algebra.

Arveson \cite{Arv67} wrote the seminal paper on this approach in 1967. He started with an ergodic, measure preserving transformation on $[0,1]$ with Lebesgue measure $m$. This induces a $*$-automorphism of $L^\infty[0,1]$. From this, he constructs a concrete weak-$*$ closed
\emph{nonselfadjoint} operator algebra that encodes the dynamics. The main result is that two such ergodic actions are \emph{conjugate} if and only if the operator algebras are unitarily equivalent. In a sequel, he and Josephson \cite{ArvJos69} recast this in the norm-closed setting of a homeomorphism on a locally compact operator space. Again they built a concrete operator algebra, and sought conditions relating conjugacy of the two dynamical systems to isomorphism of the associated algebras.

An important conceptual development was the work of Peters \cite{Pet84} in which he defined the \emph{semicrossed product} as a universal operator algebra, associated to a single $*$-endomorphism acting on an arbitrary C*-algebra. The universal nature of this definition freed the operator algebra from a specific representation. This predated the modern theory of abstract operator algebras. Peter's definition was readily generalized to the action of an arbitrary semigroup of endomorphisms of a (not necessarily selfadjoint) operator algebra. However, as we shall see, there are many choices to make when considering the class of allowable representations that do not arise when studying the action of $\bZ_+$ on a C*-algebra.

The study of nonselfadjoint operator algebras is greatly influenced by another seminal paper
of Arveson \cite{Arv69} in which he generalized Sz.Nagy's theory of dilations \cite{SFBK} to general operator algebras. The impact of this paper is widespread and profound, and we will only focus on its implications for semicrossed products. Arveson's point of view was that every operator algebra $A$ should be a subalgebra of a C*-algebra, and among all C*-algebras which can be generated by a (completely isometric) copy of $A$, there is a preferred one called the \emph{C*-envelope}. This is the analogue of the Shilov boundary of a function algebra, and is the \emph{unique minimal} C*-algebra containing $A$. The existence of this C*-envelope was established by Hamana \cite{Ham79}. A very revealing new proof was provided by Dritschel and McCullough \cite{DriMcC05}. They observed that every representation of the operator algebra can be dilated to a \emph{maximal representation}; and that these maximal representations extend uniquely to a $*$-representation of the C*-envelope. This leads to an explicit construction of the C*-envelope by finding a maximal dilation of any completely isometric representation. So to understand an operator algebra, we wish to understand the nature of the C*-envelope, and the maximal dilations.

The original work on semicrossed products dealt with an action by a single $*$-auto\-mor\-phism on an abelian C*-algebra. A few early exceptions are McAsey and Muhly \cite{McAMuh83}, who have a very general view even though they use a concrete representation, and get sharp results for semigroups that totally order the group (such as $\bR_+$); and Ling and Muhly \cite{LinMuh89} for semicrossed products over $\bZ_+^n$. As the area developed, this has expanded to the consideration of an action of a semigroup on an arbitrary operator algebra by endomorphisms. More recent examples include Donsig, Katavolos and Manoussos \cite{DKM01} and Alaimia and Peters \cite{AlaPet03} for actions of $\bZ_+^n$, the first author and Katsoulis \cite{DavKat11} and Duncan \cite{Dun08} for actions of the free semigroup $\bF_n^+$; the second author \cite{Ful12} for direct sums of positive cones in the real line; and the third author and Katsoulis \cite{KakKat12} for operator algebras associated to $n$ $*$-endomorphisms of a C*-algebra.

In this survey we focus particularly on two issues which have driven the development of semicrossed products. The first is the conjugacy problem. Can we distinguish dynamical systems (up to conjugacy) by their semicrossed products? The answer appears to be positive in many general cases. Our aim is to convince the interested reader that semicrossed products (and in general nonselfadjoint operator algebras) are not just an artifact of convenience. Indeed, even for one-variable homeomorphic classical systems, the counterexample of Hoare and Parry \cite{HoaPar66} shows that the C*-crossed product is not a complete invariant for conjugacy. Outstanding effort is required to achieve a complete encoding via C*-crossed products, and is limited to  particular systems. For example, see the deep work of Giordano, Putnam and Skau \cite{GPS95} and Giordano, Matui, Putnam, and Skau \cite{GMPS08, GMPS08, GMPS10} on Cantor systems. Moreover semicrossed products of arbitrary C*-dynamical systems over abelian semigroups are \emph{unconditionally} defined and (usually) contain a copy of the system. This is in contrast to generalized C*-crossed products, where the system is embedded by a quotient when it is non-injective, or where the existence of transfer operators is required. For example, see \cite{BRV09, Cun77, Exe03, KakKat11, Lac00, Lar10, Mur02, Pas80, Sta93}; see the work of the third author with Peters \cite{KakPet13} for a discussion.

The second issue is dilation theory and identification of the C*-envelope. In particular, we will concentrate on semigroups $P$ which are contained in a group $G$ (generated by $P$). The goal is, whenever possible, to build a C*-algebra $B$ and an action $\be$ of $G$ on $B$ by $*$-automorphisms so that the semicrossed product sits completely isometrically inside the crossed product $B\rtimes_\be G$. Then we wish to show that the C*-envelope either equals this crossed product or is a full corner of it, so that it is Morita equivalent to a crossed product. We also discuss semicrossed products by certain non-abelian semigroups. In particular we discuss Ore semigroups and free semigroups.

Connecting properties of the dynamics to simplicity of C*-crossed products is of particular interest in the theory of operator algebras (see for example \cite{Sie09}). In section~\ref{S:minimal}, we discuss how the ideal structure of the C*-envelope of a semicrossed product relates to minimality of the  underlying dynamical system.

\section{Covariant Representations and Semicrossed products}

A classical dynamical system consists of a locally compact Hausdorff space $X$ and a proper continuous map $\si$ of $X$ into itself. One is interested in various questions about how iterations of the map $\si$ evolve. The C*-algebra $\rC_0(X)$ is the obvious operator algebra that encodes $X$, and one can recover $X$ as the maximal ideal space of $\rC_0(X)$. The map $\si$ induces a $*$-endomorphism of $\rC_0(X)$ by $\al(f)(x) = f(\si(x))$. This is why we require $\si$ to be proper. Thus there is an action of $\bZ_+$ on $\rC_0(X)$ by
\[
\al_n(f)(x) = \al^n(f)(x) = f(\si^n(x)) \qfor n \ge 0,
\]
where $\si^n$ indicates the composition of $n$ copies of $\si$.

In the context of more general operator algebras, we limit our attention to completely contractive endomorphisms. When $A$ is a C*-algebra, completely contractive endomorphisms are automatically $*$-endomorphisms.

\begin{definition}
Let $P$ be a semigroup and let $A$ be an operator algebra. Let $\End(A)$ denote the completely contractive endomorphisms of $A$. A \emph{semigroup dynamical system} $(A,\alpha,P)$ consists of a semigroup homomorphism $\al\colon P \to \End(A)$.
	
When $A$ is a C*-algebra, we call $(A,\alpha, P)$ a \emph{C*-dynamical system}. When $A$ is an abelian C*-algebra, it is a \emph{classical system}.  A dynamical system is said to be \emph{unital}/\emph{injective}/\emph{surjective}/\emph{automorphic} when each $\al_s$, $s \in P$, is unital/injective/sur\-jective/automorphic.
\end{definition}

For a C*-dynamical system $(A,\al,G)$ over a group $G$, one builds the universal C*-algebra determined by the pairs $(\pi,U)$, where $\pi \colon A \to \B(H)$ is a $*$-representation, $U \colon G \rightarrow \U(H)$ is a unitary group homomorphism and
\[
\pi \alpha_g(a) = U_g \pi(a) U_g^* \qfor a\in A \AND g\in G .
\]
The crossed product $A \rtimes_\al G$ is the closure of the trigonometric polynomials on $U_g$ with coefficients from $A$. The covariance relations show that it doesn't make any difference whether the coefficients are taken on the left or on the right. Moreover this algebra is universal in the sense that whenever $\pi$ and $U$ are given satisfying the covariance relations, there is a canonical $*$-representation $\pi \times U$ of $A \rtimes_\al G$ into $\B(H)$ which restricts to $\pi$ on $A$ and to $U$ on $G$. In particular when $G$ is abelian, $\pi \times U$ is faithful if and only if $\pi$ is faithful and the pair $(\pi,U)$ admits a gauge action $\{\be_{\hat{g}} : \hat{g} \in \hat{G} \}$ over the dual group.

Given a semigroup dynamical system $(A,\alpha, P)$, the goal is to construct an operator algebra in a similar manner. The best choice of a covariant representation when $P$ is not a group is less clear because there are several equivalent formulations in the group setting which yield different results for general semigroups. The following definition follows from the original choice of Peters \cite{Pet84}, and is arguably the best choice when $P$ is abelian (see \cite{Kak11-1}). We will discuss the non-abelian case later in this section.

\begin{definition}\label{D:covariant}
Let $P$ be an abelian semigroup; and let $(A,\alpha, P)$ be a dynamical system. A \emph{covariant pair} $(\pi, T)$ for $(A,\alpha, P)$ on a Hilbert space $H$ satisfies
\begin{enumerate}
\item $\pi\colon A\rightarrow B(H)$ is a completely contractive representation of $A$;
\item $T \colon P\rightarrow B(H)$ is a representation of $P$ as contractions;
\item $\pi(a) T_s = T_s\pi\alpha_s(a)$ for all $s\in P$ and $a\in A$ (covariance relation).
\end{enumerate}
We call a covariant pair $(\pi,T)$ \emph{contractive}/\emph{isometric}/\emph{unitary} if the representation $T$ of $P$ is contractive/isometric/unitary.
\end{definition}

A semicrossed product is a universal algebra with respect to a class of covariant representations. Traditionally, semicrossed products have been defined over all contractive covariant pairs $(\pi, T)$. However, as we will see, it is often worthwhile to consider more restrictive representations of the semigroup.

\begin{definition}
Let $(A,\alpha, P)$ be a dynamical system over an abelian semigroup $P$. Let $c_{00}(P,\al,A)$ be the algebra with the linear structure of the algebraic tensor product $c_{00}(P) \otimes A$ and a multiplication given by the covariance relations:
\[
(e_s\otimes a)(e_t \otimes b) = e_{s+t}  \otimes \al_t(a)b \qfor s,t \in P \AND a,b \in A .
\]
Given a family $\F$ of covariant pairs for $(A,\alpha, P)$, define a family of seminorms on $M_n\big(c_{00}(P,\al,A)\big)$ by
\[ 
\big\| \sum_{s\in P} e_s \otimes A_s \big\| = \sup\big\{  \big\| \sum_{s\in P} (T_s \otimes I_n) \pi^{(n)}(A_s) \big\|_{\B(H^{(n)})} : (\pi,T) \in \F \big\} ,
\]
where $A_n \in M_n(A)$ and $A_n=0$ except finitely often. Let $N$ be the set of elements of $c_{00}(P,\al,A)$ with norm $0$, which is an ideal. The operator algebra completion of $c_{00}(P,\al,A)/N$ with respect to this family of seminorms yields an operator algebra $A\times_{\al}^{\F}P$ called the \emph{semicrossed product of $A$ by $P$ with respect to $\F$}.
\end{definition}

Note that even if $\F$ is not a set, the collection of norms is a set of real numbers; and the supremum is finite since 
\[ \big\| \sum_{s\in P} (T_s \otimes I_n) \pi^{(n)}(A_s) \big\|  \le \sum_{s\in P} \| A_s \| .\]
If $(\pi,T)$ is a covariant pair, then there is a representation $T\times\pi$ of $c_{00}(P) \otimes A$ into $\B(H)$ such that 
\[ (T\times\pi)(e_s \otimes a) = T_s \pi(a) \qforal a \in A \AND s \in P .\]
When $(\pi,T) \in \F$, it is clear from the definition of the seminorms that this extends to a completely contractive representation of $A\times_{\al}^{\F}P$, which is also denoted by $T\times\pi$. It is routine to select a set $\R$ of covariant pairs $(\pi,T)$ from $\F$ so that for every element of $M_n\big(c_{00}(P,\al,A)\big)$ for $n\ge1$, the supremum is obtained over the set $\R$. Then we can form the representation $\rho = \sum_{(\pi,T)\in\R}^\oplus T\times \pi$. It is easy to verify that this yields a spatial completely isometric representation with respect to this family of seminorms. So we have defined an operator algebra structure on $A\times_{\al}^{\F}P$. 

We define a few classes of covariant representations which we will consider.
\begin{definition}
(i) $A\times_\alpha P$ is the \emph{$($contractive$)$ semicrossed product} determined by all covariant pairs of $(A,\alpha, P)$;
\\\strut\quad (ii) $A\times_{\al}^{\iso}P$ is the \emph{isometric semicrossed product} determined by the isometric covariant pairs of $(A,\alpha, P)$.
\end{definition}

One can define other semicrossed products. For example, one can ask that $T$ be unitary or co-isometric as well. One can also impose different norm conditions, such as a row contractive constraint on a set of generators of $P$. One can also specify a family of representations such as the family of Fock representations which we define below.  

These choices do not always yield an algebra that contains $A$ completely isometrically.  However we normally seek representations which do contain $A$ completely isometrically. When it does not, it is generally because the endomorphisms are not faithful, but the allowable representations do not account for this. This would happen, for example, if we restricted the covariant pairs so that $T$ had to be unitary while some $\alpha_s$ has kernel. The Fock representations show that the contractive and isometric semicrossed products do always contain such a copy of $A$.

\begin{example}[Fock representations] \label{E:Fock}
Let $\pi$ be any completely contractive representation of $A$ on $H$. Let $\{e_t : t \in P\}$ be the standard orthonormal basis for $\ell^2(P)$; and define isometries on $\ell^2(P)$ by $S_s e_t = e_{s+t}$. Form a Hilbert space $\wt{H} = H \otimes \ell^2(P)$; and define a covariant pair $(\wt{\pi}, V)$ acting on $\wt{H}$ by
\[
 \wt{\pi}(a) = \diag(\pi\al_t(a))_{t \in P} \qand V_s = I \otimes S_s .
\]
Clearly, each $V_s$ is an isometry. This is  a covariance relation because
\begin{align*}
\wt\pi(a) V_s (x \otimes e_t) &= \wt\pi(a) (x \otimes e_{s+t} ) \\ &= \pi\al_{s+t}(a)x  \otimes e_{s+t} \\& = V_s (\pi\al_{s+t}(a)x  \otimes e_{t})\\ & = V_s \wt\pi\al_s(a) (x \otimes e_t) .
\end{align*}
Every operator algebra $A$ has a completely isometric representation. Therefore $A$ sits completely isometrically in  $A\times_{\al}^{\iso}P$, and thus also in  $A\times_{\al}P$. 

The commutativity of $P$ is required in these calculations. When $P$ is not abelian, one can define $S_s$ to be the right shifts $R_s e_t = e_{ts}$. Then this does yield a representation of $(A,\al,P)$ satisfying properties (i) and (iii) of a covariance relation. However this produces an antihomomorphism of $P$ into $B(\wt{H})$ because it reverses multiplication: $R_sR_t = R_{ts}$. 
\end{example}

\smallbreak
\subsection*{The semigroup $\bZ_+$}
We will write $(A,\al_1)$ instead of $(A,\al,\bZ_+)$ since the system is completely determined by $\al_1$. When $A = \rC_0(X)$ and $\al_1(f) = f\circ \si$, we write this classical system as $(X,\si)$. 

The isometric semicrossed product was defined by Peters \cite{Pet84} for C*-dynamical systems $(A,\al_1)$. Later Muhly and Solel \cite{MuhSol98} showed that every contractive covariant representation of a C*-dynamical system over $\bZ_+$ dilates to an isometric covariant representation. Therefore $A \times_\al \bZ_+ \simeq A \times_\al^{\iso} \bZ_+$ for any C*-dynamical systems $(A,\al_1)$. See  \cite{MuhSol06} for a proof avoiding the use of C*-corr\-espond\-ences. In the commutative case, Peters further showed that the Fock representation yields a (completely) isometric representation of $A \times_\al^{\iso} \bZ_+$. Peters' arguments can be adapted to the non-commutative setting (see for example \cite[Section 1]{Kak12}).

Peters introduced another important construction. Consider the direct limit $\wt{A} = \dirlim (A_n,\al_1)$ where each $A_n = A$:
\[
 \xymatrix{
 A_1 \ar[r]^{\al_1} \ar[d]_{\al_1} &A_2 \ar[r]^{\al_1} \ar[d]_{\al_1} &A_3 \ar[r]^{\al_1} \ar[d]_{\al_1} &\dots \ar[r] & \wt{A}\ar[d]^{\wt{\al}_1} \\
 A_1 \ar[r]^{\al_1}&A_2 \ar[r]^{\al_1}&A_3 \ar[r]^{\al_1}&\dots \ar[r] & \wt{A}}.
\]
There are canonical maps $\omega_n:A_n\to\wt A$, and their union is dense in $\wt A$. The direct limit of the maps $\al_1$ from $A_n$ to itself extends to an \textit{injective} $*$-endomorphism $\wt{\al}_1$ of $\wt{A}$. In addition, $\wt\al_1(\omega_n(A_n)) = \omega_{n-1}(A_{n-1})$. So $\wt{\al}_1$ has dense range, and thus is an automorphism. When $\al_1$ is injective, $(\wt{A},\wt{\al},\bZ)$ is an extension of $(A,\al,\bZ_+)$. When $\al_1$ is not injective, $(\wt{A},\wt{\al},\bZ)$ is an extension of the injective quotient system $(A/R_\al, \dot \al, \bZ_+)$ where $R_\al = \ol{\cup_n \ker\al_n}$ is the \emph{radical ideal} of this dynamical system, $q$ is the quotient of $A$ onto $A/R_\al$, and $ \dot \al$ is the induced system on the quotient.

Peters shows that, for injective C*-systems, $A \times_\al \bZ_+$ sits completely isometrically in the crossed product $\wt{A} \rtimes_{\wt \al} \bZ$. He uses this to show that every isometric covariant representation dilates to a unitary covariant representation. Thus for injective C*-dynamical systems over $\bZ_+$, the semicrossed product $A \times_\al \bZ_+$ is completely isometrically isomorphic to the semicrossed product for unitary covariant relations, $A \times_\al^{\uni} \bZ_+$. The non-injective case was studied by the third author and Katsoulis \cite{KakKat11} and the third author \cite{Kak11-1}. This requires a dilation to an injective system which we discuss  later.

The extension of these identifications to systems over nonselfadjoint operator algebras is more complicated. The third author and Katsoulis \cite{KakKat10} show that whenever $\al_1$ is a completely isometric automorphism of $A$, one gets $A \times_\al^{\iso} \bZ_+ \simeq A \times_\al^{\uni} \bZ_+$. However the contractive and the isometric semicrossed products do not coincide in general. This follows from the counterexamples of Parrott \cite{Par70} or Kaijser-Varopoulos \cite{Var74} of three commuting contractions which do not dilate to three commuting isometries. These examples yield a contactive covariant relation for $(A(\bD^2),\id,\bZ_+)$ which does not dilate to an isometric covariant relation. This problem is strongly connected to one of the fundamental problems in dilation theory: the existence of commuting lifting theorems. The reader is referred to \cite{DavKat12-2} by the first author and Katsoulis  for a full discussion.

There is a fairly general situation where a positive result holds. When $\al_1$ is a completely isometric automorphism of the tensor algebra $\T_X^+$ of a C*-corr\-espond\-ence in the sense of Muhly and Solel \cite{MuhSol98}, every contractive covariance relation dilates to an isometric one. Therefore, 
$\T_X^+ \times_\al \bZ_+ \simeq \T_X^+ \times_\al^{\iso} \bZ_+. $ 
This was established by the first author and Katsoulis \cite{DavKat10} when $\T_X^+$ is Popescu's noncommutative disc algebra $A_n$. For general C*-corr\-espond\-ences it was shown by the third author and Katsoulis \cite{KakKat10} under the assumption that $\al_1$ fixes the diagonal of $\T_X^+$. In full generality, it was proven by the first author and Katsoulis \cite[Section 12]{DavKat12-2} as a consequence of commutant lifting. An ad hoc proof of this was later added to the end of \cite{KakKat10}.

\smallbreak
\subsection*{Spanning cones}

We are concerned with what occurs when we replace $\bZ_+$ with other semigroups. A general semigroup is rather pathological; so we restrict our attention to semigroups that imbed into a group. This requires that $P$ satisfy both left and right cancellation. For an abelian semigroup, there is the well known construction of the Grothendieck group; and when the semigroup has cancellation, this imbedding is faithful. For non-abelian semigroups, it is more complicated. We discuss the case of Ore semigroups in section~\ref{S:Cenv}.

\begin{definition}
Let $G$ be an abelian group. A unital semigroup $P\subseteq G$ is called a \emph{cone}. A cone $P$ is a \emph{positive cone} if $P\cap -P = \{0\}$. A (not necessarily positive) cone $P$ is called a \emph{spanning cone} if $G=-P+P$.
\end{definition}

If $P$ is a spanning cone, we define a pre-order on $G$ by $g \leq h$ if and only if $h-g \in P$. This makes $G$ and $P$ into directed sets because if $g=-s_1+s_2$ and $h=-t_1+t_2$ in $G$ are written with $s_i,t_i\in P$, then one can see that $g$ and $h$ are both dominated by $s_2+t_2$. If in addition $P$ is a positive cone, then this yields a partial order on $G$. Conversely, given a partially ordered abelian group $G$, we obtain a positive cone $P := \{p\in G : p\geq 0\}$.

If the partial order on $G$ induced by $P$ makes $G$ into a lattice (i.e., any two elements $a,b$ in $G$ have a least upper-bound $a \vee b$ and a greatest lower-bound $a \wedge b$), then $(G,P)$ is an \emph{abelian lattice-ordered group}.

Let $(G,P)$ be an abelian lattice-ordered group. Write any $g\in G$ as $g = -s +t$ with $s,t\in P$. If we set $g_+ = t-(s \wedge t)$ and $g_-=s-(s \wedge t)$, we obtain $g_+,g_-\in P$ with $g_+\wedge g_-=0$ such that $g=-g_- + g_+$. This decomposition is unique \cite{Goo86}.

\begin{definition}
Let $(G,P)$ be a latticed-ordered abelian group. A contractive representation $T\colon P\rightarrow B(H)$ is called \emph{regular} if there exists a unitary representation $U \colon G\rightarrow B(K)$ with $H\subseteq K$ such that
\[
P_H U_g|_H=T_{g_-}^*T_{g_+} \qforal g\in G .
\]
A contractive representation $T\colon P\rightarrow B(H)$ is called \emph{Nica-covariant} if
\[
T_sT_t^*=T_t^*T_s \quad\text{whenever}\quad s,t\in P \AND s\wedge t=0.
\]
\end{definition}

The property of Nica-covariance was introduced by Nica \cite{Nic92} in the study of C*-algebras generated by isometric representations of positive cones of quasi-lattice-ordered groups. These representations have proven to be of importance in the theory of crossed product C*-algebras by semigroups, e.g., see Laca and Raeburn \cite{LacRae96}. In particular when $P = \bZ_+^n$, Nica-covariance is equivalent to saying that the canonical generators are mapped to $n$ doubly-commuting contractions, meaning that $T_iT_j^*=T_j^*T_i$ for $i \ne j$.

Let $(G,P)$ be an abelian lattice-ordered group, and let $T$ be an isometric Nica-covariant representation of $P$. Since $T$ is isometric, $Q_s:=T_sT_s^*$ is a projection for each $s\in T$. The property of Nica-covariance guarantees that
\[
Q_sQ_t=Q_{s\wedge t} \qforal s,t\in P.
\]
Thus isometric Nica-covariant representations represent $P$ not only as a semigroup, but also as a lattice.

If $U \colon G\rightarrow B(H)$ is a unitary representation of an abelian lattice-ordered group $(G,P)$, then $U|_P$ is clearly regular. Conversely, if $U$ is a unitary representation of $P$, then the unitaries automatically $*$-commute, and the representation extends to a representation of $G$ by $U_{-s+t} = U_s^*U_t$. 
In particular, any unitary representation of $P$ is Nica-covariant. 

In some instances, a Nica-covariant representation of $P$ is automatically regular. For automorphic C*-dynamical systems $(A,\alpha,\bZ_+^n)$, Ling and Muhly \cite{LinMuh89} show that Nica-covariant covariant representations are regular. The second author has generalized this result to C*-dynamical systems $(A,\alpha, P)$ where $P=\sum_{i=1}^{n\ \oplus}P_i$, with each $P_i$ a positive cone in $\bR_+$ \cite{Ful12}. In general, if $(G,P)$ is an arbitrary abelian lattice-ordered group, it is unknown whether a Nica-covariant covariant representation of a C*-dynamical system $(A,\alpha, P)$ is necessarily regular.

\begin{definition}
If $(G,P)$ is a lattice-ordered abelian group, let
$A\times_{\al}^{\nc}P$ denote the \emph{Nica-covariant semicrossed product} determined by the regular Nica-covariant covariant pairs.
\end{definition}

Every regular contractive covariant representation dilates to an isometric representation, and conversely every isometric covariant representation is regular. So there is no need for a regular semicrossed product. A regular Nica-covariant (contractive) covariant representation dilates to a Nica-covariant isometric covariant representation. When $(G,P)$ is a lattice-ordered abelian group, the Fock representation of Example~\ref{E:Fock} is Nica-covariant. Therefore $A\times_{\al}^{\nc}P$ contains $A$ completely isometrically. We show in \cite{DFK13} that the Fock representations completely norm $A\times_{\al}^{\nc}P$. 

A first attempt to generalize the one-variable results to spanning cones was carried by Duncan and Peters \cite{DunPet10} for homeomorphic classical systems. The object under consideration in \cite{DunPet10} was the concrete algebra generated by the Fock representations by point evaluation maps. We extend their results to non-commutative C*-dynamical systems in \cite{DFK13}.

Finally we mention some generalizations of And\^o's dilation of two commuting contractions. Ling and Muhly \cite{LinMuh89} show that if $(A,\alpha,\bZ_+^2)$ is a unital automorphic C*-dynamical system, then every contractive covariant representation dilates to a unitary covariant representation. In our language, this means $A\times_\al\bZ_+^2 \simeq A \times_\alpha^\uni \bZ_+^2$. Then Solel \cite{Sol06} established a general result in the setting of product systems which implies for a C*-dynamical system $(A,\alpha,\bZ_+^2)$ that $A\times_\al\bZ_+^2 \simeq A \times_\alpha^\iso \bZ_+^2$.

\smallbreak
\subsection*{Free semigroups}
Let $\bF_n^+$ be the free semigroup on $n$ generators. Since this semigroup is finitely generated, one can define variants of the semicrossed product while maintaining the action on the left as in the abelian case. These semicrossed products were first considered by the first author and Katsoulis \cite{DavKat11} for classical systems over $\bF_+^n$ and Duncan \cite{Dun08} who identifies these classical semicrossed products with a free product of simpler algebras.

\begin{definition}
Let $(A,\alpha,\bF_n^+)$ be a C*-dynamical system. A \emph{left covariant family} $(\pi, \{T_i\}_{i=1}^n)$ of $(A,\alpha,\bF_n^+)$ consists of
\begin{enumerate}
\item $\pi \colon A\rightarrow B(H)$ is a $*$-representation;
\item $T_i\in B(H)$ and $\|T_i\| \le 1$ for $1 \le i \le n$;
\item $T_i\pi(\al_i(a))=\pi(a)T_i$ for every $a\in A$ and $i=1,\ldots,n$.
\end{enumerate}
\end{definition}

The reason this does not quite fit the rubrik of Definition~\ref{D:covariant} is because items (ii) and (iii) are incompatible. The homomorphism $\al$ of $\bF_+^n$ into $\End(A)$ is determined uniquely by the family $\{\al_i\}$, namely $\al_w(a) = \al_{i_1}\al_{i_2}\dots\al_{i_k}(a)$ for $w = i_1i_2\dots i_k \in \bF_+^n$. Thus we often write this system as $(A,\{\al_i\}_{i=1}^n)$. One can also uniquely extend $\{T_i\}$ to a homomorphism $T$ of $\bF_+^n$ into $B(H)$ by $T_w = T_{i_1}T_{i_2}\dots T_{i_k}$. However the covariance relation doesn't work out as one might expect. One must use the anti-automorphism that sends $w$ to its opposite word $\ol{w}=i_k\ldots i_1$. The covariance relation can now be expressed as $\pi(a) T_{\ol{w}} = T_{\ol{w}}\, \pi\alpha_w(a)$. (Alternatively, one could save (iii) by defining $T_w$ to be $T_{i_k}T_{i_{k-1}}\dots T_{i_1}$ instead. So (iii) now holds; but the map from $w$ to $T_w$ would be an antihomomorphism, and so (ii) would fail.)

\begin{definition}
Let $(A,\{\al_i\}_{i=1}^n)$ be a C*-dynamical system.
\begin{enumerate}
\item The \emph{free semicrossed product} $A \times_\al \bF_+^n$ is the semicrossed product algebra with respect to all left covariant contractive families.
\item The tensor algebra $\T^+{(A,\al)}$ is the semicrossed product algebra with respect to the row contractive left covariant contractive families $(\pi, \{T\}_{i=1}^n)$, satisfying $\|\big[ T_1\ \dots\ T_n \big]\| \le 1$.
\end{enumerate}
\end{definition}

The tensor algebra is indeed a tensor algebra in the sense of Muhly and Solel \cite{MuhSol98} associated to a natural C*-corr\-espond\-ence. The following theorem was proved in the classical case by the first author and Katsoulis \cite[Proposition 6.2]{DavKat11}. The general case was proved by the authors in \cite{DFK13}. The dilation of each $T_i$ can be described as an explicit Schaefer type dilation.

\begin{proposition}\cite{DFK13} \label{P:dilate free}
Let $(A,\{\al_i\}_{i=1}^n)$ be a C*-dynamical system. Then every left covariant contractive family dilates to a left covariant isometric family.
\end{proposition}

In the classical case, it is also shown in \cite{DavKat11} that every row contractive covariant family dilates to a row isometric covariant family. In particular, as in Example~\ref{E:Fock}, the Fock representations using the right shifts $R_i$ on $\wt H = H \otimes \ell^2(\bF_+^n)$ are row isometric. This shows that the C*-algebra $A$ is faithfully represented inside both $A \times_\al \bF_+^n$ and the tensor algebra $\T^+{(A,\al)}$.

\smallbreak
\subsection*{Other nonabelian semigroups}
The theory of semicrossed products over non-abelian semigroups is considerably less developed. The work on the free semigroup is the main exception. As noted in the free semigroup situation, the left covariance relations lead to an antihomomorphism of $P$ into $B(H)$. This can result in the multiplication rule on $c_{00}(A,\al,P)$ failing to be associative. Another solution is needed. There is further discussion of this in \cite{DFK13}. Instead one can use the right variant of the covariance relation.  

\begin{definition}
Let $P$ be a semigroup; and let $(A,\alpha, P)$ be a dynamical system. A \emph{right covariant representation} of $(A,\alpha, P)$ on a Hilbert space $H$ is a pair $(\pi, T)$ such that
\begin{enumerate}
\item $\pi\colon A\rightarrow B(H)$ is a completely contractive representation of $A$;
\item $T \colon P\rightarrow B(H)$ is a representation of $P$ as contractions;
\item $T_s \pi(a) = \pi\alpha_s(a)T_s$ for all $s\in P$ and $a\in A$ (right covariance relation).
\end{enumerate}
We call a covariant pair $(\pi,T)$ \emph{contractive}/\emph{isometric}/\emph{unitary} if the representation $T$ of $P$ is contractive/isometric/unitary.
\end{definition}

A semicrossed product can then be obtained in a similar manner.

\begin{definition}
Let $(A,\alpha, P)$ be a dynamical system over a non-abelian semigroup $P$. Define $c_{00}(P,\al,A)_r$ to be an algebra with the  linear structure of the algebraic tensor product $A \otimes c_{00}(P)$ and a multiplication given by the right covariance relations:
\[
(a \otimes e_s)(b \otimes e_t) = a\al_s(b) \otimes e_{st} \qfor s,t \in P \AND a,b \in A .
\]
Given a family $\F$ of right covariant representations, define a seminorm on the matrix algebras $M_n\big( c_{00}(P,\al,A)_r \big)$ by
\[
 \big\| \sum_{s\in P} A_s \otimes e_s \big\| = \sup \big\{ \big\| \sum_{s\in P} \pi^{(n)}(A_s) (I_n \otimes T_s) \big\|_{\B(H^{(n)})} : (\pi,T) \in \F  \big\} ,
\]
where $A_s \in M_n(A)$ and $A_s = 0$ except finitely often. The operator algebra completion $\fA(A,\al,\F)_r$ of $c_{00}(P,\al,A)_r$ with respect to this family of seminorms is the \emph{semicrossed product of $A$ by $P$ with respect to $\F$}. 

We write $\fA(A,\al,P)_r$ when $\F$ is the family of all covariance relations, and $\fA(A,\al,\iso)_r$ when $\F$ is the family of isometric covariance relations.
\end{definition}

We write elements of $P$ on the right to correspond with the right covariance relations. As before, the supremum is over a set of real numbers, even when the family $\F$ is not a set; and this supremum is finite because it is dominated by $\sum_{s\in P} \|A_s\|$. Each covariant pair $(\pi,T)$ in $\F$ yields a completely contractive representation of $\fA(A,\al,\F)_r$ determined by $(\pi\times T)(a \otimes e_s) = \pi(a) T_s$. Again one can select a set $\R$ from $\F$ so that the supremum of the norms is attained over $\R$ for every element. The direct sum of the representations $\pi\times T$ over $(\pi,T)\in\R$ yields a completely isometric isomorphism of $\fA(A,\al,\F)_r$ into some $B(H)$. 

In the non-abelian case, even being a subsemigroup of a group is not sufficient structure for strong dilation results. We will show better results can be obtained for Ore semigroups.

\section{Semicrossed products and the conjugacy problem}

We will now discuss how semicrossed products can be used to distinguish between non-conjugate dynamical systems.  It is a classical problem to determine invariants that can discriminate between two inequivalent dynamical systems. Historically this is the first application of semicrossed products.
 
\smallbreak
\subsection*{The semigroup $\bZ_+$.}

First recall what it means for two classical dynamical systems to be  conjugate. 
\begin{definition}
Two classical dynamical systems $(X, \si)$ and $(Y, \tau)$ are said to be \emph{conjugate} provided that there is a homeomorphism $\gamma$ of $X$ onto $Y$ such that $\tau\gamma = \gamma\si$.
\end{definition}

In \cite{HoaPar66}, Hoare and Parry give an example of a homeomorphic classical system $(X,\si)$ such that $\si$ and $\si^{-1}$ are not conjugate. However, the crossed product C*-algebras $C(X) \rtimes_\si \bZ$ and $C(X) \rtimes_{\si^{-1}}\bZ$ are always isomorphic. Thus, in general, C*-crossed products do not provide a complete isomorphism invariant for conjugacy of dynamical systems.

In contrast, semicrossed products have proven to very useful for this purpose. Arveson, in his seminal paper \cite{Arv67}, constructed a concrete weak-$*$ closed analogue of the semicrossed product to distinguish ergodic transformations. Then with Josephson \cite{ArvJos69}, he tackled the norm-closed case and obtained important partial results.

\begin{example}[Arveson-Josephson algebras]\label{ex: arv-jos}
Let $X$ be a locally compact, Hausdorff space and let $\phi \colon X\rightarrow X$ be a homeomorphism. Assume that there is a separable regular Borel probability measure $m$ on $X$ such that
\begin{enumerate}
\item $m\circ\phi<<m$ (quasi-invariance);
\item $m(U)>0$ for all non-empty open subsets $U$ of $X$;
\item the set of periodic points $P=\bigcup_{n>0}\{\phi^{n}(x)=x\}$ has measure zero.
\end{enumerate}
It is noted in \cite{ArvJos69} that the existence of such a measure $m$ is not automatic.

Let $U$ be the unitary on $L^2(X,m)$ defined by
\[	
Uf=\big(\tfrac{d\,m\circ\phi}{dm}\big)^{1/2} f\circ\phi .
\]
For $f\in C_0(X)$, let $L_f$ be the multiplication operator $L_f h=fh$ on $L^2(X,m)$.	
The Arveson-Josephson algebra $\mathfrak{A}(X,\phi, m)$ is the norm-closed algebra generated by all finite sums
\[
L_{f_0}+L_{f_1}U+\ldots+L_{f_n}U^n \qfor f_i\in C_0(X) .
\]
Note that $U$ satisfies the covariance relation $UL_f U^*=L_{f\circ\phi}$.  Peters \cite[Example IV.10]{Pet84} showed that the Arveson-Josephson algebra is isomorphic to the semicrossed product $C(X_0)\times_\phi\bZ_+$.
\end{example}

Arveson and Josephson answered the conjugacy problem under the assumption that the measure is ergodic and invariant under the homeomorphism $\phi$.

\begin{theorem}[Arveson-Josephson]
Let $\mathfrak{A}(X_1,\phi_1, m_1)$ and $\mathfrak{A}(X_2,\phi_2, m_2)$ be Arv\-eson-Josephson algebras. Further suppose that $m_2$ is ergodic and invariant under $\phi_2$. Then the following are equivalent
\begin{enumerate}
\item  $(X_1,\phi_1)$ and $(X_2,\phi_2)$ are conjugate.
\item $\mathfrak{A}(X_1,\phi_1, m_1)$ and $\mathfrak{A}(X_2,\phi_2, m_2)$ are isometrically isomorphic.
\item $\mathfrak{A}(X_1,\phi_1, m_1)$ and $\mathfrak{A}(X_2,\phi_2, m_2)$ are algebraically isomorphic.
\end{enumerate}	
\end{theorem}

Peters \cite{Pet84} improved on this by allowing $\phi$ to be a continuous map of $X$ into itself, but $X$ had to be compact and there could be no periodic points at all. Hadwin and Hoover \cite{HadHoo88}, still in the compact case, weakened the condition on fixed points considerably. In 2008, the first author and Katsoulis \cite{DavKat08} removed all such conditions on the dynamical systems and proved the following.

\begin{theorem}[Davidson-Katsoulis]
Let $X_i$ be a locally compact Hausdorff spaces, and let $\sigma_i$ be a proper continuous map of $X_i$ into itself. Then the following are equivalent:
\begin{enumerate}
\item $(X_1,\sigma_1)$ and $(X_2,\sigma_2)$ are conjugate.
\item $C_0(X_1)\times_{\sigma_1}\bZ_+$ and $C_0(X_2)\times_{\sigma_2}\bZ_+$ are completely isometrically isomorphic.
\item $C_0(X_1)\times_{\sigma_1}\bZ_+$ and $C_0(X_2)\times_{\sigma_2}\bZ_+$ are algebraically isomorphic.	
\end{enumerate}
\end{theorem}

An example of Buske and Peters \cite{BusPet98} shows that if $\si$ is a conformal map of the unit disc onto itself with a unique fixed point interior to the disc, and $\al(f)(z)= f(\si(z))$ is the induced automorphism of the disc algebra $A(\bD)$, then $A(\bD) \times_\al \bZ_+ \simeq A(\bD) \times_{\al^{-1}} \bZ_+$. The results of \cite{DavKat08} also deal with the isomorphism problem for semicrossed products of certain function algebras in the plane. There is an analogous result to the one above, with the pathology of the Buske-Peters example being the only obstruction to conjugacy.

Now we turn to non-abelian C*-dynamical systems over $\bZ_+$. The appropriate notion of equivalence is the following.

\begin{definition}
Two C*-dynamical systems $(A,\al_1)$ and $(B,\be_1)$ are said to be \emph{outer conjugate} if there is a $*$-isomorphism $\ga$ of $A$ onto $B$ and a unitary $v$ in the multiplier algebra $\M(A)$ such that
\[
\al_1(a) = \ad_v \ga^{-1}\be_1\ga(a) = v(\ga^{-1}\be_1\ga(a))v^* \qforal a\in A .
\]
\end{definition}

Suppose that $(A,\al_1)$ and $(B,\be_1)$ are outer conjugate. Let $\Bs$ and $\Bt$ be the generators of $\bZ_+$ in $A\times_\al \bZ_+$ and $B\times_\be \bZ_+$, respectively. As $\ga$ is a $*$-isomorphism of $A$ onto $B$, it extends to a $*$-isomorphism of their multiplier algebras. Thus we can define $w = \ga(v^*)$. A computation shows that the map sending $a$ to $\ga(a)$ and $\Bs a$ to $\Bt w \ga(a)$ extends to a completely isometric isomorphism of the semicrossed products $A\times_\al \bZ_+$ and $B\times_\be \bZ_+$. We are interested in when the converse is valid.

Semicrossed products often provide a complete isometric isomorphic invariant for outer conjugacy of C*-dynamical systems. Muhly and Solel \cite[Theorem 4.1]{MuhSol00} show that if $A=B$ is a separable C*-algebra, $\al_1$ and $\be_1$ are automorphisms, and one of them has full Connes spectrum, then the semicrossed product are isometrically isomorphic if and only if the systems are outer conjugate. The first author and Katsoulis \cite{DavKat08-2} establish this result for separable, simple C*-algebras, also when $\al_1$ is an automorphism. The methods in both of these papers are more difficult than the following much more general results due to the first and third authors \cite{DavKak12}.

\begin{theorem}[Davidson-Kakariadis]\label{T: DavKak conj}
Let $(A,\al_1)$ and $(B,\be_1)$ be unital C*-dynam\-ical systems.
Suppose that any one of the following holds:
\begin{itemize}
\item $\al_1$ is injective,
\item $\al_1$ is surjective,
\item $A$ has trivial centre, $($e.g., when $A$ is simple$)$,
\item $A$ is abelian,
\item $\al_1(A)' \cap A$  is finite $($no proper isometries$)$.
\end{itemize}
Then $A\times_\al \bZ_+$ and $B\times_\be \bZ_+$ are isometrically isomorphic
if and only if $(A,\al_1)$ and $(B,\be_1)$ are outer conjugate.
\end{theorem}

\smallbreak
\subsection*{Multivariable systems}

A multivariable classical dynamical system is a locally compact, Hausdorff space $X$ together with a collection $\sigma_1,\dots,\sigma_n$ of proper continuous maps of $X$ into itself. For simplicity we write $(X,\si)$ for $(X, \{\si_i\}_{i=1}^n)$. We say that $(X,\si)$ has multiplicity $n_\si := n$. Such systems were studied by the first author and Katsoulis \cite{DavKat11}.

There are no labels on the maps $\si_i$. So conjugacy of the two systems should allow for an arbitrary permutation of the maps by an element of the symmetric group $S_n$. It is less obvious, but possible in some circumstances, that one can switch continuously from one permutation to another. This leads to the following definition.

\begin{definition}[Davidson-Katsoulis]
Two multivariable dynamical systems $(X,\sigma)$ and $(Y,\tau)$ are \textit{piecewise conjugate} if there is a homeomorphism $\gamma\colon X \to Y$ and an open cover $\{ \U_\alpha : \alpha \in S_n \}$ of $X$ so that
\[
\tau_i \gamma|_{\U_\alpha} = \gamma \sigma_{\alpha(i)}|_{\U_\alpha}  \qfor \alpha \in S_n .
\]
\end{definition}

The following example illustrates how this can occur.
\begin{example}
Let $X=Y=[0,1]$ and $n=2$. Let 
\begin{alignat*}{2}
\si_1(x) &= 
\begin{cases}
 2x &\ \FOR 0 \le x \le \tfrac13\\
 1-x &\ \FOR \tfrac13 \le x \le  1
\end{cases}
&\qand
\si_2(x) &= 
\begin{cases}
 1-x &\ \FOR 0 \le x \le \tfrac23\\
 2x-1 &\ \FOR \tfrac23 \le x \le 1
\end{cases}
\\
\tau_1(x) &= 
\begin{cases}
 2x &\ \FOR 0 \le x \le \tfrac13\\
 1-x &\ \FOR \tfrac13 \le x \le \tfrac23\\
 2x-1 &\ \FOR \tfrac23 \le x \le 1
\end{cases}
&\qand
\tau_2(x) &= \ \ 
1-x \qquad  \FOR 0 \le x \le 1 . 
\end{alignat*}
Then let $\U_\id = [0,\tfrac23)$ and $\U_{(12)} = (\tfrac13,1]$. Observe that $\si_1$ and $\tau_1$ agree on $\U_{\id}$, as do $\si_2$ and $\tau_2$, while $\si_1$ and $\tau_2$ agree on $\U_{(12)}$ as do $\si_2$ and $\tau_1$. The two systems $([0,1],\{ \si_1,\si_2\})$ and $([0,1],\{\tau_1,\tau_2\})$ are piecewise conjugate, but not conjugate.
\end{example}

The semigroup that naturally acts here is the free semigroup $\bF_+^n$. The covariance relations can be expressed in terms of the generators of $\bF_+^n$ by
\[  \pi(f) T_i = T_i \pi(f\circ \si_i) \qforal f \in \rC_0(X) \AND 1 \le i \le n .\]
Two semicrossed products are considered in that paper, the free semicrossed product $\rC_0(X) \times_\sigma \bF_+^n$ and the tensor algebra $\T^+{(X,\si)}$. It appears that piecewise conjugacy is strongly connected to the structure of these algebras. The main result of \cite{DavKat11} is the following.

\begin{theorem}[Davidson-Katsoulis]\label{T: DK pc}
Let $(X,\sigma)$ and $(Y,\tau)$ be two multivariable dynamical systems. If $\rC_0(X) \times_\sigma \bF_+^n$ and $\rC_0(Y) \times_\tau \bF_+^n$ or the tensor algebras $\T^+{(X,\si)}$ and $\T^+{(Y,\tau)}$ are isomorphic, then $(X,\sigma)$ and $(Y,\tau)$  are piecewise conjugate.
\end{theorem}

In the case $n=1$, the converse is easy. The converse remains an open question for the semicrossed product for $n\ge2$. For the tensor algebra, the converse has been verified in special cases such as $\Dim X \le 1$ or $n\le 4$  \cite{DavKat11, Ram09}. The tensor algebra seems to be easier to work with than other algebras defined for this kind of system.

This theory of piecewise conjugacy has proven to be useful outside of operator algebras. Cornelissen \cite{Cor13}, and Cornelissen and Marcolli \cite{CorMar11, CorMar13} apply the notion of piecewise conjugacy and Theorem \ref{T: DK pc} to obtain results in number theory and graph theory. In particular, in \cite{CorMar11} they obtain a group isomorphism invariant for their systems and in \cite{CorMar13} they classify Banach algebras associated to quantum mechanical systems, as in Bost-Connes-Marcolli systems \cite{BosCon95, HaPau05}, applied to the reconstruction of graphs. It becomes apparent that a converse of Theorem \ref{T: DK pc} is of importance, and it is established for certain systems that they study.

The problem of the existence of a complete invariant for multivariable systems was studied in the non-commutative context by the third author and Katsoulis \cite{KakKat12}. They define a notion of piecewise conjugacy for automorphic C*-dynamical systems in terms of the induced action $(\hat{A},\hat{\al})$ on the Fell spectra $\hat{A}$ (see \cite{Ern75} for a useful characterization). A non-commutative version of Theorem \ref{T: DK pc} is then obtained. Their examination included other algebras, not all of which are semicrossed products. They include the free semicrossed product, the tensor algebra, the Nica-covariant semicrossed product $A \times_\al^{\nc} \bZ_+^n$ and the isometric semicrossed product $A \times_\al^{\iso} \bZ_+^n$, as well as examples formed by using as a prototype the operator algebras of~\cite{DRS11} related to analytic varieties. The following result applies to this general class of algebras, but we only state it for the two semicrossed products that we have been discussing.

\begin{theorem}[Kakariadis-Katsoulis]\label{T: KakKat pc}
Let $(A,\al)$ and $(B,\be)$ be automorphic unital multivariable C*-dynamical systems. If either the two semicrossed products $A \times_\al \bF_+^n$ and $B \times_\be \bF_+^n$ or the two tensor algebras $\T^+(A,\al)$ and $\T^+(B,\be)$ are isometrically isomorphic, then the multivariable systems $(\hat{A},\hat{\al})$ and $(\hat{B},\hat{\be})$  are piecewise conjugate.
\end{theorem}

In this generality, a converse does not hold even for $n=1$. Kadison and Ringrose \cite{KadRin67} show that there exists a homogeneous C*-algebra $A$ and an automorphism $\al$ of  $A$ which is universally weakly inner but not inner. If the converse of Theorem \ref{T: KakKat pc} were valid for tensor algebras, then $A \times_{\al} \bZ^+$ and $A\times_\id \bZ^+$ would be isomorphic and hence would be outer conjugate by Theorem \ref{T: DavKak conj}. But this would imply that $\al$ is inner, a contradiction.

There are other alternatives for extending outer conjugacy from the one variable case to the multivariable level. Recall that every multivariable system gives rise to a C*-corr\-espond\-ence $X{(A,\al)}$ \cite{MuhSol98}.

\begin{definition}
Let $(A,\al)$ and $(B,\be)$ be multivariable C*-dynamical systems. We say that $X{(A,\al)}$ and $X{(B,\be)}$ are \emph{unitarily equivalent} when there is a $*$-isomorphism $\ga\colon A \rightarrow B$ and a $n_\be \times n_\al$ unitary matrix $[u_{ij}]$ with coefficients in $\M(B)$ such that
\[
\diag\{\be_i\ga\}_{ i=1}^{n_\be} \cdot [u_{ij}] = [u_{ij}] \cdot \diag\{\ga \al_j\}_{ j=1}^{ n_\al}.
\]
When $[u_{ij}]$ is diagonal up to a permutation, and thus $n_\al=n_\be$, the multivariable systems $(A,\al)$ and $(B,\be)$ are said to be \emph{outer conjugate}.
\end{definition}

Again the first part of the following theorem is valid for a much wider collection of algebras associated to the dynamical system.

\begin{theorem}[Kakariadis-Katsoulis]\label{T: KakKat st}
Let $(A,\al)$ and $(B,\be)$ be unital multivariable C*-dynamical systems. Suppose that either $A$ is stably finite or 
 $\al_i$ and $\be_j$ are $*$-epimorphisms for all $i,j$. Then
\begin{enumerate}
\item If $A \times_\al \bF_+^n$ and $B \times_\be \bF_+^n$ are isometrically isomorphic, then $n_\al=n_\be$ and the correspondences $X{(A,\al)}$ and $X{(B,\be)}$ are unitarily equivalent.
\item $\T^+(A,\al)$ and $\T^+(B,\be)$ are isometrically isomorphic if and only if the correspondences $X{(A,\al)}$ and $X{(B,\be)}$ are unitarily equivalent.
\end{enumerate}
\end{theorem}

This result implies that the tensor algebra is a complete isometric isomorphic invariant for unitary equivalence of classical systems. Thus one way to show that the tensor algebras are a complete isomorphic invariant for piecewise conjugacy would be to identify the two equivalence relations. In particular we have the following corollary.

\begin{corollary}[Kakariadis-Katsoulis]
Let $(A,\al)$ and $(B,\be)$ be two automorphic multivariable C*-dynamical systems and assume that $A$ has trivial center. Then the following are equivalent:
\begin{enumerate}
\item $A\times_{\al} \bF_{n_{\al}}^{+}$  and $B\times_{\be} \bF_{n_{\be}}^{+}$ are isometrically isomorphic.
\item $\T^+{(A,\al)}$ and $\T^+{(B,\be)}$ are isometrically isomorphic.
\item $X{(A,\al)}$ and $X{(B,\be)}$ are unitarily equivalent.
\item $(A,\al)$ and $(B,\be)$ are outer conjugate.
\end{enumerate}
\end{corollary}

One consequence of the results in \cite{DavKat11} and \cite{KakKat12} is that isomorphism of the tensor algebras implies equality of the multiplicities of the systems. The following example \cite{KakKat12} shows that the converse does not hold. However the conjugacy problem remains open: are the tensor algebras a complete isometric isomorphic invariant for unitary equivalence in general?

\begin{example}\label{E:counter}
Let $A=B=\O_2= \ca(s_1,s_2)$ be the Cuntz algebra on two generators, let $\al=(\al_1, \al_2)$, with $\al_1=\al_2=\id$, and let $\be(x) = s_1xs_1^* + s_2xs_2^*$, for all $x \in \O_2$. Then $\begin{bmatrix} s_1 & s_2 \end{bmatrix}$ is a unitary matrix in $M_{1,2}(\O_2)$ and
\[
\big[ \be(x) \big] \begin{bmatrix} s_1 & s_2 \end{bmatrix}
= \begin{bmatrix} s_1 & s_2 \end{bmatrix} \begin{bmatrix} \al_1(x) & 0 \\ 0 & \al_2(x) \end{bmatrix}
\qforal x\in \O_2 .
\]
Hence the tensor algebras $\T^{+} (A, \al)$ and $\T^{+} (B ,\be) $ are completely isometrically isomorphic. However the multiplicities $n_\al$ and $n_\be$ are not the same.
\end{example}

\section{C*-envelopes of Semicrossed Products} \label{S:Cenv}

Arveson \cite{Arv69} developed dilation theory to study the representations of an arbitrary operator algebra. Today, this can be interpreted to mean an abstract operator algebra in the sense of Blecher, Ruan and Sinclair \cite{BRS}. At the time, it meant a subalgebra of a C*-algebra with the inherited matrix norm structure. He observed that an operator algebra can sit (completely isometrically) inside a variety of C*-algebras. In the case of a function algebra $A$, restriction to the \v{S}ilov boundary $Y$ yields the minimal C*-algebra $\rC(Y)$ containing $A$ isometrically. The C*-envelope was conceived as the non-commutative analogue.

Let $j\colon A\rightarrow C$ be a completely isometric isomorphism of $A$ into a C*-algebra $C$ such that $C=\ca(j(A))$. Then $C$ is a \emph{C*-cover} of $A$. The C*-envelope of $A$ is a  C*-cover  $\cenv(A)=\ca(\iota(A))$ of $A$ with the property that whenever $C=\ca(j(A))$ is any other C*-cover of $A$,  there exists a surjective $*$-homomorphism $\Phi \colon C\to \cenv(A)$ such that $\Phi \circ j=\iota$.

Whilst Arveson calculated the C*-envelope of many operator algebras \cite{Arv69,Arv72}, these early papers did not establish the existence of the C*-envelope in general. He was particularly interested in irreducible representations $\rho$ of $\ca(A)$ with the property that $\rho|_A$ has a unique completely positive extension to $\ca(A)$. These are called \emph{boundary representations}, and they are the analogue of points in the Choquet boundary of a function algebra. Existence of the C*-envelope was established by Hamana \cite{Ham79}. Hamana's proof was extremely significant, but his methods did not connect the C*-envelope to dilation theory.

In 2001, Dritschel and McCullough \cite{DriMcC05} announced an exciting new proof of the C*-envelope that shed new light on this issue. They called a (completely contractive) representation $\rho$ of $A$  \emph{maximal}  if all dilations of $\rho$ have the form $\si = \rho \oplus \si'$. They established that every representation does dilate to a maximal one, and that the maximal representations have a unique completely positive extension to a $*$-representation of $\ca(A)$ that factors through the C*-envelope. Thus the C*-envelope can be obtained by taking any completely isometric representation of $A$ and dilating it to a maximal representation $\rho$; then $\cenv(A) = \ca(\rho(A))$.

Their paper did not address the issue of \emph{irreducible} representations with this property. Arveson \cite{Arv08} revisited this issue in light of these new developments, and succeeded in establishing the existence of sufficiently many boundary representations in the separable case. Arveson's vision was finally fully realized when the first author and Kennedy \cite{DavKen13} showed that every operator system has a noncommutative Choquet boundary.

Identifying the C*-envelope of an operator algebra has been a central problem in dilation theory for the last 45 years. The work of Dritschel and McCullough \cite{DriMcC05} focussed attention on the notion of a maximal dilation. This motivated the third author and Peters \cite{KakPet13} to outline a program on the C*-envelope. For semicrossed products this leads to the issue of dilating non-invertible systems to automorphic systems, and connecting the C*-envelope of a semicrossed product to a C*-crossed product. Such results imply that the C*-envelope of a semicrossed product is a reasonable choice for generalizing C*-crossed product constructions to non-invertible dynamical systems.

\smallbreak
\subsection*{C*-envelopes for systems over $\bZ_+$}

The first computation of the C*-envelope of a semicrossed product was obtained by Muhly and Solel \cite[Corollary 6.9]{MuhSol98}. They were interested in semicrossed products over $\bZ_+$ as a special case of their study of tensor algebras of C*-corr\-espond\-ences. They established (in somewhat different language) a prototype of the kind of result we are seeking. They noted that Peters \cite{Pet84} was close to this result in 1984.

\begin{theorem}[Muhly-Solel] \label{MS_C*_auto}
Let $A$ be a unital C*-algebra and let $\al$ be a $*$-automorphism. Then
\[	
\cenv(A\times_\al\bZ_+) \simeq A\rtimes_\al\bZ .
\]
\end{theorem}

Peters \cite{Pet08} approached the topic again in 2008, where he studied the C*-env\-el\-opes of semicrossed products of classical dynamical systems over $\bZ_+$. His motivation in this work was to explore the relationship between C*-envelopes of semicrossed products and crossed product C*-algebras. Again, a strong relationship was found. Let $X$ be a compact Hausdorff space and let $\sigma$ be a continuous, surjective map on $X$. From the dynamical system $(X,\sigma,\bZ_+)$, one can construct a homeomorphic dynamical system $(\wt{X},\wt{\sigma},\bZ)$ as a projective limit:
\[
\xymatrix{	X & X \ar[l]_{\sigma} & X \ar[l]_{\sigma}& X \ar[l]_{\sigma} &  \ldots \ar[l]_{\sigma} & \wt{X} \ar[l]}
\]
Viewing $X$ as a subset of $\wt{X}$ and $\wt{\si}$ as a homeomorphic extension of $\si$, Peters shows that $C(X)\times_\sigma\bZ_+$ embeds completely isometrically into $C(\wt{X})\rtimes_{\wt{\si}}\bZ$. This leads us to the following result.

\begin{theorem}[Peters]
Let $(X, \si, \bZ_+)$ be a dynamical system where $\phi$ is continuous and surjective. Then
\[	
\cenv(C(X)\times_\si\bZ_+) \simeq C(\wt{X})\rtimes_{\wt{\si}}\bZ ,
\]
where $(\wt{X}, \wt{\si},\bZ)$ is the system described above.
\end{theorem}

The system $(\wt{X}, \wt{\si},\bZ)$ is no other than the minimal automorphic direct limit extension $(X,\si,\bZ_+)$. It was observed by the third author and Katsoulis \cite{KakKat10} that Peters' result can be extended to non-classical systems. Soon after in \cite{ KakKat11}, they treated the general case in the following way. Define the ideal
\[
I = (\ker\al)^\perp = \{ a \in A : a \ker\al = \{0\} \} ;
\]
and let $q\colon A \to A/I$ be the quotient. Let $B = A \oplus (A/I\otimes c_0)$, and write a typical element of $B$ as
\[
 b = a_0 \otimes e_0 + \sum_{n\ge1} q(a_n) \otimes e_n
 \qfor  a_n \in A.
\]
Define a $*$-endomorphism $\be$ on $B$ by
\[
 \be \big(a_0 \otimes e_0 + \sum_{n\ge1} q(a_n) \otimes e_n \big) =
 \al(a_0)  \otimes e_0 +  \sum_{n\ge 1} q(a_{n-1}) \otimes e_n .
\]
It is not hard to check that this map is injective and dilates $\al$. Therefore $(B,\be,\bZ_+)$ is injective and admits a minimal automorphic extension $(\wt B, \wt\be, \bZ_+)$.

\begin{theorem}[Kakariadis-Katsoulis]
Let $(A,\al, \bZ_+)$ be a C*-dynamical system, and let $(B,\be, \bZ_+)$ and  $(\wt B, \wt\be, \bZ)$ be constructed as above. Then the C*-envelope $\cenv(A \times_\al \bZ_+)$ is a full corner of $\wt B \rtimes_{\wt\be} \bZ$.
\end{theorem}

This result is expressed in the language of C*-corr\-espond\-ences in \cite{KakKat11}, and the third author provides an ad hoc proof in \cite{Kak11-1}. The system $(B,\be,\bZ_+)$ is obtained by ``adding a tail'' which is reminiscent of the theory of graph algebras. This idea was first introduced by Muhly and Tomforde \cite{MuhTom04} to study non-injective C*-corr\-espond\-ences. However their construction doesn't always yield an injective C*-dynamical system \cite[Proposition 3.14]{KakKat11}. The third author and Katsoulis \cite{KakKat11} show that one can choose from a variety of tails in order to preserve the appropriate structure of the system.

The semicrossed product of a C*-dynamical system over $\bZ_+$ is an example of a C*-corr\-espond\-ence in the sense of Muhly and Solel \cite{MuhSol98}. Therefore one can get the C*-envelope via the Cuntz-Pimsner algebra. This was accomplished under certain conditions by Fowler, Muhly and Raeburn \cite{FowMuhRae03}; and was extended by Katsoulis and Kribs \cite{KatKri06} to show that the C*-envelope of the tensor algebra $\T_X^+$ of a C*-corr\-espond\-ence is the Cuntz-Pimsner algebra in the sense of Katsura \cite{Kat04}. Therefore the C*-envelope of $A \times_\al \bZ_+$ is the universal C*-algebra generated by isometric covariant pairs $(\pi,V)$ such that
\[
\pi(a)(I-VV^*)=0 \qfor a \in (\ker\al )^\perp .
\]
However, this characterization does not provide the association to a C*-crossed product that we are seeking.

The examination of similar results for systems of nonselfadjoint operator algebras is trickier, because of the counterexamples of Parrott \cite{Par70} and Kaijser-Varopoulos \cite{Var74}. Since 2008, papers by the first author and Katsoulis \cite{DavKat10, DavKat12-1, DavKat12-2} and the third author and Katsoulis \cite{KakKat10} have focused on the cases where one can calculate the C*-envelope of a semicrossed product over $\bZ_+$. These results rely on an old observation of Arveson \cite{Arv69} that every completely isometric automorphism of an operator algebra $A$ extends to a $*$-automorphism of $\cenv(A)$.

\begin{theorem}[Kakariadis-Katsoulis] 
Let $A$ be an operator algebra and let $\alpha$ be a completely isometric automorphism. Then
\[
\cenv(A\times_\al^\iso \bZ_+)\simeq \cenv(A)\rtimes_\alpha \bZ.
\]
More generally, if $\alpha$ is a completely contractive endomorphism of $A$ which extends to a $*$-endomorphism of $\cenv(A)$, then
\[
\cenv(A\times_\al^\iso \bZ_+)\simeq \cenv(\cenv(A)\times_\alpha\bZ_+).
\]
\end{theorem}

In particular when $\al$ is a completely isometric automorphism of a tensor algebra $\T_X^+$, then $\cenv(\T_X^+ \times_\al \bZ_+) \simeq \O_X \rtimes_\al \bZ$.

\smallbreak
\subsection*{C*-envelopes for dynamical systems over abelian semigroups}

In contrast to the one-variable case, one cannot generally obtain such unconditional results for other semigroups. This fails even for trivial systems over $\bZ_+^n$ since one cannot dilate three commuting contractions to unitaries \cite{Par70, Var74}.  Ling and Muhly's automorphic And\^{o} theorem \cite{LinMuh89}, mentioned earlier, can be used to show that if $(A,\alpha,\bZ_+^2)$ is a unital automorphic C*-dynamical system, then $\cenv(A\times_\al\bZ_+^2)\simeq A\rtimes_\alpha \bZ^2$.

One of the first attempts to compute the C*-envelope of a semicrossed product over more general positive cones is due to the second author in \cite{Ful12}. The semigroups studied there were of the form $P=\sum_{i=1}^{\oplus n} P_i$, where each $P_i$ was a positive cone in $\bR_+$. If $G$ is the subgroup of $\bR^n$ generated by $P$, then $(G,P)$ is a lattice-ordered abelian group. Thus it is natural to consider the Nica-covariant representations of this semigroup. In this instance Nica-covariant representations of $P$ are automatically regular. The methods in \cite{Ful12} generalize those developed in \cite{KakKat10} to the multivariable context.

\begin{theorem}[Fuller]
Let $P=\sum_{i=1}^{n \oplus}P_i$, where each $P_i$ is a positive cone in $\bR_+$. If $(A,\alpha, P)$ is a dynamical system consisting of completely isometric automorphisms of an operator algebra $A$, then
\[
\cenv(A\times_\alpha^{\nc} P) \simeq \cenv(A)\rtimes_\alpha G.
\]
\end{theorem}

For spanning cones of more general abelian groups without a lattice structure, we have the following.

\begin{theorem}\cite{DFK13}
Let $P$ be a spanning cone of an abelian group $G$. If $(A,\alpha, P)$ is a dynamical system consisting of completely isometric automorphisms, then
\[
\cenv(A\times_\alpha^{\iso} P) \simeq \cenv(A)\rtimes_\alpha G .
\]
\end{theorem}

For the rest of this section, we will restrict our attention to C*-dynamical systems. Peters' method \cite{Pet08} of creating a homeomorphic extension of a classical system extend to C*-dynamical systems over spanning cones. This method has been developed by Laca \cite{Lac00}.

Since spanning cones are directed, we can construct the direct limit C*-algebra $\wt{A}$ associated to the connecting $*$-homomorphisms $\al_t \colon A_s \to A_{s +t}$ for $A_s :=A$ and $s, t\in P$. For every $p \in P$ a $*$-automorphism $\wt{\al}_p$ is defined with respect to the diagram
\[
\xymatrix{
A_s \ar[r]^{\al_t} \ar[d]^{\al_p} & A_{s+t} \ar[d]^{\al_p} \ar[r] &\wt{A} \ar[d]^{\wt\al_p} \\
A_s \ar[r]^{\al_t} & A_{s+t}  \ar[r] &\wt{A}
}
\]
Furthermore the mapping $\wt{\al} \colon P \to \Aut(\wt{A})$ is a semigroup homomorphism. Therefore it extends to a group homomorphism of $G= - P +P$. Thus we obtain an automorphic system $(\widetilde{A}, \widetilde{\al}, G)$, which is called the \emph{minimal automorphic extension} of $(A,\al, P)$.

When $(A,\al,P)$ is injective, the system $(\wt{A},\wt{\al},G)$ is indeed an extension. However, as mentioned earlier, the $*$-homomorphisms of $A_s$ into $\wt{A}$ factor through the quotient by the radical ideal $R_\alpha = \overline{\bigcup_{s\in P}\ker\al_s}$. Define $\dot{A}=A/R_\alpha$. As $R_\alpha$ is invariant under each $\alpha_s$, this induces an injective dynamical system $(\dot{A},\dot{\alpha},P)$.

We will concentrate on the case when $(G,P)$ is an abelian lattice-ordered group. Because the Nica-covariant representations take into account the lattice structure of the semigroup, we consider $A\times_\alpha^{\nc} P$ as the natural semicrossed product in this context.

\begin{theorem}\cite{DFK13}
Let $(G,P)$ be a lattice-ordered abelian group, and let $(A,\alpha, P)$ be an injective C*-dynamical system. Let $(\widetilde{A}, \widetilde{\al}, G)$ be its minimal automorphic extension. Then
\[
\cenv(A\times_\alpha^{\nc} P)\simeq \widetilde{A}\rtimes_{\widetilde{\al}}G.
\]
\end{theorem}

We can say more when $P=\bZ_+^n$. Note that, in this case, all Nica-covariant representations are automatically regular. We will be able to remove the condition that the $*$-endomorphisms $\alpha_s$ are injective. In this case, the C*-envelope turns out to be a full corner of the crossed product. The first step is an \emph{adding tails} method to imbed $(A,\al,\bZ_+^n)$ into an injective system. It is much more complicated than the $\bZ_+$ case. There are many ways to extend the system to an injective one, but we need to construct a minimal dilation that will provide the C*-envelope. 

We will write $\bo1, \dots, \Bi, \dots, \Bn$ for the standard generators of $\bZ_+^n$.  Also write $\un{0} = (0,\dots,0)$ and $\un 1 = (1,\dots,1)$. For $\un x = (x_1,\dots,x_n) \in \bZ_+^n$, define
\[
 \supp(\un x) = \{\Bi : x_i\ne 0 \}  \qand
  \un x^\perp = \{\un y \in \bZ_+^n : \supp(\un y) \cap \supp(\un x) = \mt \} .
\]
For each $\un x \in \bZ_+^n\setminus\{\un{0}\}$, consider the ideal $\big( \bigcap_{\Bi\in \supp(\un x)} \ker\al_\Bi \big)^\perp$. We require the largest ideal contained in this ideal which is invariant for $\al_{\un y}$ whenever $\un y \in \un x^\perp$, namely
\[
 I_{\un x} = \bigcap_{\un y \in \un x^\perp} \al_{\un y}^{-1}
 \Big(\big( \bigcap_{\Bi\in \supp(\un x)} \ker\al_\Bi \big)^\perp\Big) .
\]
In particular, $I_\un{0} = \{0\}$, and $I_{\un x} = I_{\un 1} = \big( \bigcap_{\Bi=\bo1}^\Bn \ker\al_\Bi \big)^\perp$ for all $\un x \ge \un 1$.

Define $B_{\un x} := A/I_{\un x}$ for $\un x \in \bZ_+^n$ and let $q_\un x$ be the quotient map of $A$ onto $B_\un x$. Set $B = \sum_{\un x \in \bZ_+^n}\!\!\oplus\, B_\un x$. A typical element of $B$ will be denoted as
\[ b = \sum_{\un x\in\bZ_+^n} q_{\un x}(a_{\un x}) \otimes e_{\un x} ,\]
where $a_{\un x} \in A$. Observe that $I_{\un x}$ is invariant under $\al_\Bi$ when $x_i=0$. Therefore we can define $*$-endomorphisms $\be_{\Bi} \in \End(B)$ as follows:
\begin{align*}
 \be_{\Bi}(q_{\un x}(a) \otimes e_{\un x} ) =
 \begin{cases}
 q_{\un x}\al_\Bi(a) \otimes e_{\un x} + q_{\un x + \Bi}(a) \otimes e_{\un x + \Bi} & \text{ if } x_i=0,\\
 q_{\un x}(a) \otimes e_{\un x + \Bi} & \text{ for } x_i\ge 1.
 \end{cases}
\end{align*}
It is clear that the compression of $(B,\be,\bZ_+^n)$ to $A$ is the system $(A,\al,\bZ_+^n)$. So $(B,\be,\bZ_+^n)$ is a dilation of $(A,\al,\bZ_+^n)$. Furthermore $(B,\be,\bZ_+^n)$ is injective, hence it admits a minimal automorphic extension $(\wt B, \wt\be, \bZ^n)$.

\begin{theorem}\label{NC Zn}\cite{DFK13}
Let $(A,\alpha,\bZ_+^n)$ be a non-injective C*-dynamical system. Let $(B,\beta,\bZ_+^n)$ be the injective dilation of $(A,\al,\bZ_+^n)$ described above, and let $(\widetilde{B},\widetilde{\beta}, \bZ^n)$ be the automorphic extension of $(B,\beta,\bZ_+^n)$. Then $\cenv(A\times_{\al}^{\nc}\bZ_+^n)$ is a full corner of $\widetilde{B}\rtimes_{\widetilde{\beta}}\bZ^n$.
\end{theorem}

As a consequence of our methods we also obtain the following characterization for the C*-envelope.

\begin{theorem}\cite{DFK13}
Let $(A,\alpha,\bZ_+^n)$ be a  C*-dynamical system over $\bZ_+^n$. Then the C*-envelope of $A \times_\al^{\nc} \bZ_+^n$ is the universal C*-algebra generated by Nica-covariant isometric pairs $(\pi,V)$ such that
\[
\pi(a) \cdot \prod_{\Bi \in \supp(\un x)} (I - V_\Bi V_\Bi^*) = 0 \qforal a \in I_{\un x}.
\]
\end{theorem}

Combining this with the works of Sims and Yeend \cite{SimYee10}, and Carlsen, Larsen, Sims and Vittadello \cite{CLSV11}, and with a gauge invariant uniqueness theorem we provide in \cite{DFK13}, we get that the C*-envelope is in fact the Cuntz-Nica-Pimsner algebra of a relative product system.

\smallbreak
\subsection*{Ore Semigroups}
One class of non-abelian semigroups with a very useful structure are the Ore semigroups. They imbed into a group in a special way. Dilations of actions over Ore semigroups have been studied by Laca \cite{Lac00} rather successfully.

\begin{definition}
A cancellative semigroup $P$ is an \emph{Ore semigroup} if $Ps\cap Pt\neq\emptyset$ for all $s,t\in P$.
\end{definition}

It was shown by Ore \cite{Ore31} and Dubreil \cite{Dub43} that any Ore semigroup $P$ can be embedded into a group $G$ such that $G=P^{-1}P$. Laca \cite{Lac00} shows that $G$ is characterized by the following universal property: if $\varphi \colon P \rightarrow K$ is a semigroup homomorphism into a group $K$, then it extends to a (necessarily unique) group homomorphism of $G$ into $K$.

Ore semigroups come with an order structure. For $s,t$ in $P$, we say $s\leq t$ if $t\in Ps$. This ordering directs $P$, as the Ore property ensures that given $s_1,s_2$ in $P$,
there is a $t$ in $P$ such that $s_1 \le t$ and $s_2\le t$. Therefore a direct limit $\wt{A}$ can be constructed with respect to the connecting maps $\al_t \colon A_s \to A_{ts}$ for $A_s:=A$ and $s,t$ in $P$. Then an automorphic C*-dynamical system $(\wt{A},\wt{\al},G)$ can be obtained by carefully defining the endomorphisms $\wt{\al}_p$. Let $w_s \colon A_s \to \wt{A}$ and define $\wt{\al}_p$ by the rule
\[
\wt{\al}_p w_s(a) = w_q \al_t(a) \quad\text{ when } ts=qp, \FORAL a\in A.
\]
When $\al$ is not injective, $(A,\al,P)$ imbeds into $(\wt{A},\wt{\al},G)$ via a quotient by the radical ideal $R_\al = \ol{\cup_s \ker\al_s}$.

This construction was first obtained by Laca \cite{Lac00} for injective systems. A similar construction shows also that an isometric semigroup action $V \colon P \rightarrow \B(H)$ of an Ore semigroup $P$ dilates to a unitary group action $U \colon G \rightarrow \B(K)$ of the universal group $G$ of $P$. We obtain the following characterization of the C*-envelope.

\begin{theorem} \label{thm: Ore} \cite{DFK13}
Let $(A,\alpha,P)$ be a unital C*-dynamical system over an Ore semigroup.
Then
\[
\cenv(\fA(A,P,\iso)_r) \simeq \widetilde{A}\rtimes_{\widetilde{\alpha}}G ,
\]
where $(\widetilde{A},\widetilde{\alpha},G)$ is the minimal automorphic extension of $(\dot A, \dot\al, P)$.
\end{theorem}

The condition that the system is unital is crucial here as it enables us to use the identity of $A$ as the identity of the semicrossed product. Hence it is generated by a copy of $A$ and a copy of $P$.

In the non-unital case we can use a folklore idea to realize this by restricting ourselves to \emph{non-degenerate} representations of $A$. That is, we define $\fA(A_{\nd},P,\iso)_r$ to be the universal algebra for isometric right covariant pairs $(\pi,V)$  such that $\pi$ is in addition a non-degenerate representation of $A$.

\begin{theorem} \cite{DFK13}
Let $A$ be a unital C*-algebra, and let $(A,\alpha,P)$ be a $($possibly non-unital$)$ C*-dynamical system over an Ore semigroup $P$. Let $(\wt{A},\wt{\alpha},G)$ be the minimal automorphic extension of $(A,\al,P)$. Then $\cenv(\fA(A_{\nd},P,\iso)_r)$ is a full corner of $\widetilde{A}\rtimes_{\widetilde{\alpha}}G$.
\end{theorem}

\smallbreak
\subsection*{Free semigroups}
Let $(A,\alpha,\bF_n^+)$ be a C*-dynamical system over $\bF_n^+$. The calculation of the C*-envelope of $A\times_\al\bF_n^+$ depends on the dilation Proposition~\ref{P:dilate free}. We have the following description of the C*-envelope in the automorphic case. The key step in the proof is to first dilate the $T_i$ to isometries, and then further extend them to unitaries while maintaining the covariance relations.

\begin{theorem}
Let $(A,\{\al_i\}_{i=1}^n)$ be a unital automorphic C*-dynamical system. Then
\[
\cenv(A\times_\al \bF_+^n) \simeq A \rtimes_\al \bF^n.
\]
\end{theorem}

The C*-envelope of a tensor algebra is shown to be the Cuntz-Pimsner algebra of the relative C*-corr\-espond\-ence by Katsoulis and Kribs \cite{KatKri06}, extending Fowler, Muhly and Raeburn  \cite{FowMuhRae03}, who established this for the case of a faithful and strict correspondence. It can also be associated to a generalized crossed product that resembles those considered by Paschke \cite{Pas80}. This procedure for classical systems was carried out by the first author and Roydor \cite{DavRoy11}. It was extended to the non-commutative case by the third author and Katsoulis \cite{KakKat11}. Most of these results follow by applying a version of the ``adding tails'' technique to imbed non-injective systems into injective ones.

\section{Minimality and C*-envelopes of semicrossed products} \label{S:minimal}

Earlier, we argued that the C*-envelope of a semicrossed product is the appropriate analogue of a generalized group C*-crossed product. This is partly due to its universal property. It is also due in part to the fact that, much like C*-crossed products, minimality of a dynamical system can sometimes be detected by the C*-envelope of the semicrossed product. This is a consequence of the dilation theory, which connects properties of the group action to properties of the semigroup action.

\begin{definition}
A C*-dynamical system $(A,\alpha, P)$ is \emph{minimal} if $A$ does not contain any non-trivial $\alpha$-invariant ideals.
\end{definition}

When $P$ is a group $G$,  this notion of minimality coincides with the usual one for automorphic systems. Furthermore, when $(G,P)$ is a lattice-ordered abelian group, a unital automorphic C*-dynamical system $(A,\alpha,P)$ is minimal if and only if $(A,\alpha, G)$ is minimal.

We recall some definitions from the theory of crossed products. Let $\al \colon G \rightarrow \Aut(A)$ be a group homomorphism of an abelian group $G$. The universal property of the crossed product implies there is a \emph{gauge action} of the dual group $\wh{G}$ by $*$-automorphisms which acts on finite sums by
\[
\ga_{\hat g}( \sum_{g \in G} U_g a_g ) = \sum_{g \in G} \ip{g,\hat g} U_g a_g .
\]
Consider the conditional expectation $E \colon A \rtimes_\al G \rightarrow A$ given by
\[
 E(F) = \int_{\wh{G}} \ga_{\hat{g}}(F) d\hat{g},
\]
and define the Fourier coefficients of an element $F$ in the crossed product by $E_g(F) :=  E(U_{-g} F)$. 

Moreover, if $\I$ is an $\ad_{U^*}$-invariant ideal of $A \rtimes_\al G$ over $G$, then $E_{g}(\I)$ is an $\al$-invariant ideal of $A$ over $G$. An ideal $\I$ of $A \rtimes_\al G$ is called \emph{Fourier-invariant} if $E_{g}(\I) \subseteq \I$ for all $g\in G$.

\begin{theorem}\cite{DFK13}
Let $(G,P)$ be a lattice-ordered abelian group and let $(A,\alpha, P)$ be a unital injective C*-dynamical system. Then the following are equivalent:
\begin{enumerate}
\item $(A,\alpha, G)$ is minimal;
\item $(\widetilde{A},\widetilde{\alpha},G)$ is minimal;
\item $\cenv(A\times_\alpha^{\nc}P)\simeq \widetilde{A}\rtimes_{\widetilde{\alpha}} G$ has no non-trivial Fourier-invariant ideals.
\end{enumerate}
\end{theorem}

With the help of Theorem \ref{NC Zn}, one can remove the hypothesis of injectivity above in the case when $P=\bZ_+^n$. The conclusion is the same because minimality of any of $(A,\alpha, \bZ_+^n)$ , $(B,\beta,\bZ_+^n)$ or $(\widetilde{B},\widetilde{\beta},\bZ^n)$ implies injectivity of the system.

Hence, if $\cenv(A\times_\alpha^{\nc}\bZ_+^n)$ is simple, then $(A,\alpha,\bZ_n^+)$ is minimal. The converse however is false even for $n=1$. Indeed let $A$ be any simple C*-algebra; then $(A,\id_A,\bZ_+)$ is trivially minimal but $A \rtimes_\al \bZ = A \otimes \rC(\bT)$ is not simple. Such systems fail to be topological free, which explains in part why minimality does not imply simplicity. The classical dynamical systems are more amenable, and stronger results can be obtained.

\begin{definition}
A classical dynamical system $(X,\varphi,P)$ over a semigroup $P$ is called \emph{topologically free} if $\{x\in X : \varphi_s(x) \neq \vpi_r(x)\}$ is dense in $X$ for all $s,r \in P$.
\end{definition}

This definition is a reformulation of group topological freeness appropriate for the semigroup context. The following result is obtained by observing that topological freeness is equivalent to injectivity of the semigroup action, followed by an application of a result of Archbold and Spielberg \cite{ArcSpi93}. We will let $(\wt{X},\wt{\varphi},G)$ denote the minimal automorphic extension of $(X,\varphi,P)$. Note that surjectivity of the maps $\varphi_s$ is equivalent to the injectivity of the corresponding endomorphisms.

\begin{theorem}\cite{DFK13}\label{T: min P}
Let $(X,\phi,P)$ be a surjective classical system over a lattice-ordered abelian group $(G,P)$. Then the following are equivalent:
\begin{enumerate}
\item $(X, \varphi, P)$ is minimal and $\varphi_s \neq \varphi_r$ for all $s, r \in P$;
\item $(\wt{X}, \wt{\varphi}, G)$ is minimal and topologically free;
\item the C*-envelope $\rm{C}(\wt{X}) \rtimes_{\wt{\varphi}} G$ of $\rm{C}(X) \times_\varphi^{\nc} P$ is simple.
\end{enumerate}
\end{theorem}

If $(X,\phi,\bZ_+^n)$ is a classical system, then we can drop the hypothesis of injectivity by applying Theorem \ref{NC Zn}.  Again, we can conclude that simplicity of $\rm{C}(\wt{Y}) \rtimes_{\wt{\tau}} \bZ^n$ is equivalent to the surjective system $(\wt{Y}, \wt{\tau}, \bZ^n)$ being minimal and topologically free, and it forces $(X, \varphi, \bZ_+^n)$ to be surjective.


\end{document}